\documentclass[a4paper,10pt]{amsart}
\usepackage{latexsym,amssymb,amsmath,amsthm,amsfonts}
\usepackage[T1]{fontenc}
\usepackage[latin1]{inputenc}
\usepackage{bbm}
\usepackage{graphicx}
\usepackage{color} 
\numberwithin{equation}{section}
\usepackage{caption}

\newtheorem{thm}{Theorem}[section]
\newtheorem{prop}[thm]{Proposition}
\newtheorem{lemma}[thm]{Lemma}
\newtheorem{cor}[thm]{Corollary}
\theoremstyle{definition}
\newtheorem{defi}[thm]{Definition}

\newcommand{\dem}{\noindent \textbf{Proof: }}

\newcommand{\findem}{\vspace{-.55cm} \begin{flushright} $\square~$
\end{flushright} \vspace{.2cm} }

\newcommand{\merge}{\operatorname{merge}}
\newcommand{\spli}{\operatorname{split}}

\newcommand{\lb}{[\![}
\newcommand{\rb}{]\!]}

\def \sous#1#2{\mathrel{\mathop{\kern 0pt#1}\limits_{#2}}}

\newcommand{\mA}{\mathcal{A}}

\newcommand{\al}{\alpha}
\newcommand{\be}{\beta}
\newcommand{\ga}{\gamma}
\newcommand{\de}{\delta}
\newcommand{\la}{\lambda}

\title{The structure of alternative tableaux.}
\author{Philippe Nadeau}
\address{Faculty of Mathematics, University of Vienna, Nordbergstra{\ss}e 15, 1090 Vienna,
Austria.\\
E-mail: {\tt philippe.nadeau@univie.ac.at}.}

\begin{document}

 \begin{abstract}
 In this paper we study alternative tableaux introduced by Viennot \cite{VienCamb}. These tableaux are in simple bijection with permutation tableaux, defined previously by Postnikov \cite{Postnikov}. We exhibit a simple recursive structure for alternative tableaux. From this decomposition, we can easily deduce a number of enumerative results. We also give bijections between these tableaux and certain classes of labeled trees. Finally, we exhibit a bijection with permutations, and relate it to some other bijections that already appeared in the literature.
 \end{abstract}

\maketitle

\section*{Introduction}

Alternative tableaux are certain fillings of Ferrers diagrams introduced by Xavier Viennot \cite{VienCamb}, in simple bijection with permutation tableaux introduced by Postnikov \cite{Postnikov} in his study of the totally positive part of the grassmannian. These permutation tableaux have then been considered by several authors.

Alternative tableaux are related to the stationary distribution of a certain Markov process from statistical physics, the Asymmetric Simple Exclusion Process (ASEP). They are more precisely connected to a {\em Matrix Ansatz} describing a general solution for this stationary distribution (see Proposition \ref{prop:connection}). Formulas for this distribution have been first computed by Sasamoto, Uchiyama and Wadati \cite{Uchi}, and involve the famous Askey-Wilson orthogonal polynomials. An understanding of alternative tableaux thus seem to be a key in order to get a better grasp of these polynomials, for which no combinatorial interpretation is known. The articles \cite{CW,CW1,CW2} develop such connections between tableaux and some special cases of the ASEP model. 

 Another line of research is related to permutations; the starting point here is that permutation tableaux of size $n$ are counted simply by $n!$. They were first studied from a combinatorial point of view in \cite{SW}, in which a bijection between these tableaux and permutations of $\{1,\ldots,n\}$ was studied in detail. Since then, other bijections have been described \cite{Bu,CN,V}. These tableaux also showed up in the work of Lam and Williams \cite{LW}  where permutation tableaux were shown to fit naturally into the {\em type A} case of the classification of Coxeter systems.
\medskip

In this article,  we show that alternative tableaux admit a natural recursive structure, which is best expressed when considering alternative tableaux as labeled combinatorial objects. The central part of this work is thus Section~\ref{sect:struct}, in which  we exhibit this recursive decomposition. These structural results are then applied in the following sections, first by giving enumerative results in a very straightforward manner, and then by encoding the recursive decomposition by certain classes of labeled trees, which are then themselves in bijection with permutations.
\medskip

Let us give a more precise outline of the paper. We introduce some elementary definitions and properties concerning alternative and permutation tableaux in Section~\ref{sect:first}. We then give our main results relative to the structure of alternative tableaux in Section~\ref{sect:struct}, the recursive structure being a consequence of Proposition~\ref{prop:cut} and Theorem~\ref{th:decomp} in particular. Using this decomposition, we prove several enumeration results in Section~\ref{sect:enum}: this gives in particular elementary proofs of certain results of \cite{CN,CW2,LW}, as well as some new results.  We then describe how the recursive decomposition is naturally associated with certain labeled trees in Section~\ref{sect:alttrees}. Finally we exhibit in Section~\ref{sect:perm} a bijection from alternative tableaux to permutations, and stress its connection to bijections which have already appeared in the literature.

\section{Tableaux}
\label{sect:first}

\subsection{Shapes and tableaux}

\begin{figure}[!ht]
 \centering
 \includegraphics[height=4cm]{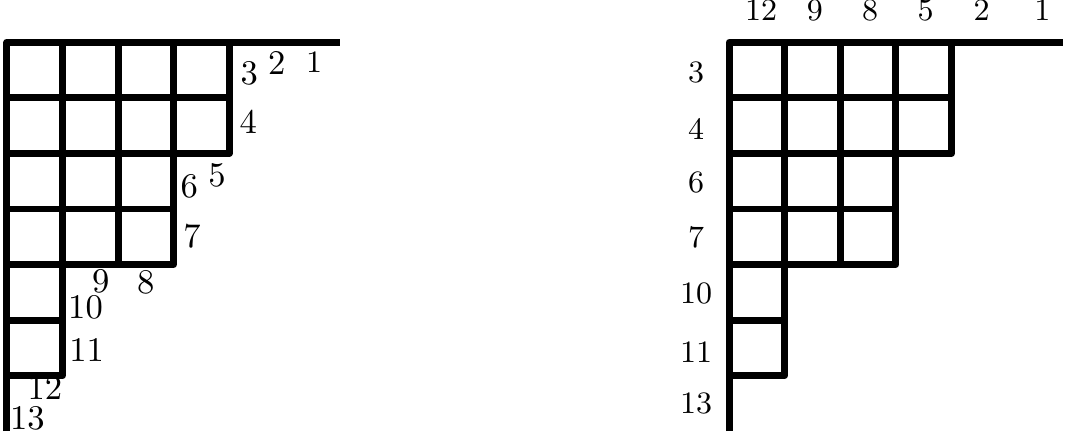}
 \caption{A shape with its standard labeling.
 \label{fig:shape}}
\end{figure}

We call \textit{shape} a staircase diagram (also called Ferrers shape) with possible empty rows or columns, cf. Figure~\ref{fig:shape}. The {\em length} of a shape is the number of rows plus the number of columns of the shape. Note that a shape is determined by its south east border, which is the path from the top right corner of the shape to its bottom left corner; it is the path labeled by the integers $1,2,\ldots 13$ on the left of Figure~\ref{fig:shape}. There are thus $2^n$ shapes of length $n$, since one can choose to go down or left at each step.

Rows and columns of shapes will be labeled by integers in the following manner: let $S$ be a shape of length $n$, and $L=\{i_1<i_2<\ldots<i_n\}$  a set of integers. Then $S$ is \textit{labeled by $L$} if the numbers $i_\ell$ are attached to the rows and columns, in increasing order following the south east border of $S$ from top right to bottom left. If $L=\{1,\ldots,n\}$, then we say that $S$ has the \textit{standard labeling}; this is the case of the shape of Figure~\ref{fig:shape}. Note that although the labeling is always defined with respect to the south east border, we will actually write the labels on the top and left side for an improved readability, as shown on the right of Figure~\ref{fig:shape}.

Given a labeled shape, we will sometimes say row $i$ or column $j$ when we actually refer to the row with the label $i$ or the column with the label $j$. Then the cell lying at the intersection of row $i$ and column $j$ is denoted by $(i,j)$, where we have $i<j$ necessarily by definition of the labelings.
\medskip

We now define two possible ways to fill these shapes, called permutation tableaux and alternative tableaux: the two are in bijection by Viennot's Theorem~\ref{th:permalt}. What we will show in this paper is that it is better to study alternative tableaux when it comes to discover the intrinsic structure of these combinatorial objects.

\begin{defi}[Permutation tableau]
A  {\em permutation tableau} $T$ is a shape with a filling of each of its cells by $0$ or $1$ such that the following properties hold: (i) each column contains at least one $1$, and (ii) there is no cell filled by a $0$ which has simultaneously a $1$ above it in the same column and a $1$ to its left in the same row.
\end{defi}

\begin{figure}[!ht]
 \centering
 \includegraphics[height=5cm]{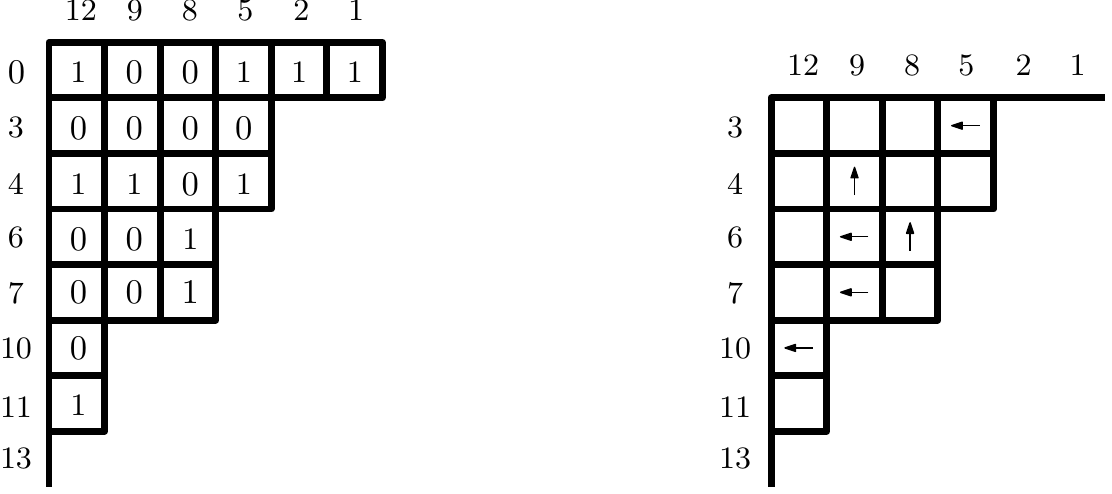}
 \caption{Bijection between permutation tableaux and alternative tableaux.
 \label{fig:permtabaltab}}
\end{figure}

Note that permutation tableaux cannot have any empty columns because of the first condition. A $0$ in a permutation tableau is  {\it restricted} if there is a $1$ above it in the same column; it is {\em rightmost restricted} if it is the rightmost such $0$ in its row. A row is {\it unrestricted} if it does not contain a restricted $0$.  A $1$ is {\it superfluous} if, somewhere above it in the same column, there is another cell containing a $1$. A permutation tableau is represented on the left of Figure~\ref{fig:permtabaltab}; lines $0,4,11$ and $13$ are unrestricted, the cells in the top row filled by $1$ appear in columns $1,2,5$ and $12$, and there are four superfluous ones in cells  $(4,5),(4,12)$ and $(11,12)$.

\begin{defi}[Alternative tableau]
\label{def:altab}
An {\em alternative tableau} is a shape with a partial filling
of the cells with left arrows $\leftarrow$ and up arrows $\uparrow$, such that all cells left of a left arrow, or above an up arrow are empty. In other words, all cells pointed by an arrow must be empty.
\end{defi}

In an alternative tableau, a {\em free row}  is a row with no left arrow, and a {\em free column} is a column with no up arrow. Thus rows (respectively columns) that are not free are in bijection with left (resp. up) arrows. A {\em free cell} is a cell which is not filled, and such that there exists no left arrow to its right and no up arrow under it; in other words, the cell is empty and no arrow points toward it. We will let $frow(T),fcol(T)$ and $fcell(T)$ denote the number of free rows, free columns and free cells of a given tableau. For the tableau $T_0$ which is represented on the right of Figure~\ref{fig:permtabaltab}, the free rows are $4,11$ and $13$ while the free columns are $1,2,5$ and $12$. There are four free cells, namely $(4,5),(4,12),(7,8)$ and $(11,12)$. Thus we have $frow(T_0)=3$ and $fcol(T_0)=fcell(T_0)=4$.

We can now state the fundamental result of Xavier Viennot, showing that alternative tableaux are actually a new simple encoding of permutation tableaux:

\begin{thm}[\cite{VienCamb}]
\label{th:permalt}
There is a bijection $\alpha$ between permutation tableaux of length $n+1$
and alternative tableaux of length $n$. If $P$ is labeled by $L$, and
 $L'$ is the set $L$ minus its smallest element, then we label $\alpha(P)$
 by $L'$, and we have:
 \begin{itemize}
 \item columns of $P$ with a $1$ in their top row correspond
 to free columns of $\alpha(P)$;
 \item unrestricted rows of $P$ (the top one excepted) correspond to free rows of $\alpha(P)$,
\item and cells of $P$ filled with superfluous $1$ correspond to free cells of $\alpha(P)$.
\end{itemize}
\end{thm}

\dem We just give a description of the bijection and its inverse, and refer to \cite{VienCamb} for more details about the proof. Given a permutation tableau $P$, transform all non superfluous $1$ to up arrows, and all rightmost restricted $0$ to left arrows; then erase all the remaining $0$ and $1$, and finally remove the first row from $P$; the result is $\alpha(P)$. An illustration is given on Figure~\ref{fig:permtabaltab}.

For the inverse bijection, given an alternative tableau $T$, add a new top row on top of it, and fill by $1$ all cells of this row that lie above a free column of $T$; then change all up arrows and free cells to $1$, and all remaining cells to $0$. The resulting tableau is $\alpha^{(-1)}(T)$.
\findem

\subsection{Alternative tableaux and the ASEP}

An important application of permutation tableaux, due to Corteel and Williams in a series of papers \cite{CW,CW1,CW2}, is related to a certain model of statistical mechanics, the ASEP: we will briefly talk about some of the connections.
The ASEP model(Asymmetric Simple Exclusion Process) is a model that can be described as the following Markov chain (see \cite{DS}). Let $\alpha,\beta,\gamma,\delta,q$ be real numbers in $[0,1]$, and $n$ a nonnegative integer. The states of the Markov chain are the $2^n$ words of length $n$ on the symbols $\circ$ and $\bullet$. Positions in the words represent sites, which can be either empty ($\circ$) or occupied by a particle ($\bullet$). The transition probabilities $p(s_1,s_2)$ between two states $s_1$ and $s_2$ model the way particles can jump from site to site, enter or exit the system:
\begin{itemize}
\item If $s_1=A\bullet\circ B$ and $s_2=A\circ\bullet B$, then $p(s_1,s_2)=\frac{1}{n+1}$ and $p(s_2,s_1)=\frac{q}{n+1}$.
\item If $s_1=A\bullet$ and $s_2=A\circ$, then $p(s_1,s_2)=\beta$ and $p(s_2,s_1)=\delta$.
\item If $s_1=\circ B$ and $s_2=\bullet B$, then $p(s_1,s_2)=\alpha$ and
$p(s_2,s_1)=\gamma$.
\item If $s_1\neq s_2$ do not correspond to any of these cases, then we set $p(s_1,s_2)=0$.
\item Finally, we have naturally $p(s_1,s_1)=1-\sum_{s_2\neq s_1}. p(s_1,s_2)$.
\end{itemize}

\begin{figure}
\centering
\includegraphics{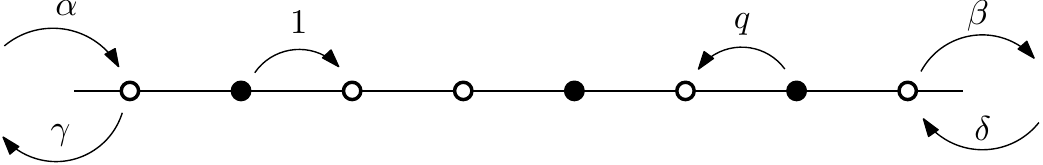}
\caption{\label{fig:asep} Illustration of the ASEP model}
\end{figure}

 The model is {\em simple} because there can be at most one particle in each site. We illustrate schematically this model in Figure~\ref{fig:asep}.

It was shown by Derrida \emph{et al.}\cite{Derrida1} that this model has a unique stationary distribution, and that moreover this distribution could be computed through the following {\em Matrix Ansatz}: suppose that we can find two matrices $D,E$, a column vector $V$ and a row vector $W$ such that the following relations hold:
\begin{equation}
\label{eq:matrel}
\begin{cases}
DE=qED+D+E\\
(\be D-\de E)V= V \\
W(\al E-\ga D)=W
\end{cases}
\end{equation}

Now let $s$ be a state of the ASEP, and let $s(E,D)$ be the word of length $n$ in $D$ and $E$ obtained through the substitutions $\bullet\mapsto D$, $ \circ\mapsto E$: for instance, the state $\circ\bullet\bullet\circ\bullet$ is associated to the word $EDDED$. Then Derrida \emph{et al} show :

\begin{prop}
\label{prop:proba}
If $D,E,V,W$ satisfy the relations~\eqref{eq:matrel}, and words in $D,E$ are interpreted as matrix products, the  probability $P_n(s)$ to be in state $s$ in  is given by
\begin{equation}
\label{eq:matrixansatz}
P_n(s)= \frac{Ws(E,D)V}{Z_n}~~~\text{with}~~~Z_n=W(D+E)^nV
\end{equation}

\end{prop}

If we have a word $w$ in the letters $E$ and $D$, we can create a shape by reading the word from left to right and interpreting each $D$ as a south step and each E as an east step, thus defining the south east boundary of a shape $\la(w)$; for instance the shape of Figure~\ref{fig:shape} is associated to the word $EEDDEDDEEDDED$; if $s$ is a state of the ASEP, we will write simply $\la(s)$ for the shape $\la(s(E,D))$. Now we can state the connection with alternative tableaux, first noticed by Corteel and Williams in a series of papers \cite{CW,CW1} and expressed in terms of permutation tableaux, and later reformulated by Viennot \cite{VienCamb} in the following way:

\begin{prop}
\label{prop:connection}
If $D,E$ are matrices that verify the first relation in~\eqref{eq:matrel}, and $w=w(E,D)$ is any word in $E,D$, then we have the following identity:

$$w=\sum_T q^{fcell(T)}E^{fcol(T)}D^{frow(T)}$$

where the sum is over all alternative tableaux of shape $\la(w)$.
\end{prop}

In \cite{Derrida1} it was in fact shown that in the case $\ga=\de=0$ there exists matrices $D,E$ verifying~\eqref{eq:matrel}. But note that in this particular case, the vectors $W$ and $V$ become respectively left and right eigenvectors for $D$ and $E$. So from Propositions~\ref{prop:proba} and~\ref{prop:connection} we obtain the following:

\begin{cor}
When  $\ga=\de=0$ we have:
\begin{equation*}
P_n(s)=\frac{\sum_{T\text{~of shape~}\la(s)} q^{fcell(T)}\al^{-fcol(T)}\be^{-frow(T)}}{\sum_{T\text{~of size~}n} q^{fcell(T)}\al^{-fcol(T)}\be^{-frow(T)}}
\end{equation*}
\end{cor}

In the ASEP model where $\ga,\de$ are general but $q=1$, Corteel and Williams found a similar expression for the stationary probabilities in terms of certain enriched alternative tableaux, see Corollary 4.2 in \cite{CW2}, by slightly generalizing the Matrix Ansatz. 

\section{The structure of alternative tableaux}
\label{sect:struct}

We define different operations on tableaux, and use them to exhibit a natural recursive structure on alternative tableaux.

\subsection{First properties of alternative tableaux}
\label{sub:firstprop}

 We denote by $\mA(n)$ the set of alternative tableaux of length $n$,
 and $\mA_{i,j}(n)$ those with $i$ free rows and $j$ free columns.
  We also denote by $\mA_{i,*}(n)$ and $\mA_{*,j}(n)$ the tableaux
  having $i$ free rows and $j$ free columns respectively.

\subsubsection{Transposition}
\label{subsub:tr}
We let $tr$ be the operation of transposing a tableau, which is
the reflection across the {\em main diagonal}, i.e. the line going south east from the top left corner; in this reflection, we naturally exchange up and left arrows. We have the following immediate result, which we state as a proposition for future reference:

 \begin{prop}
 \label{prop:transp}
  Transposition is an involution on alternative tableaux. For all $n,i,j\geq 0$,it exchanges  $\mA_{i,j}(n)$ and $\mA_{j,i}(n)$.
 \end{prop}

In fact it is easily checked that the transposition operation coincides with the involution $I$ defined in Section 7 of \cite{CW1} for permutation tableaux: more precisely, if $P$ is a permutation tableau, then we have $tr\circ\alpha(P)=\alpha\circ I(P)$, where $\alpha$ is the bijection between permutation tableaux and alternative tableaux of Theorem~\ref{th:permalt}. Note then that the trivial result on alternative tableaux from Proposition~\ref{prop:transp} demanded a much greater effort in \cite{CW1} where the authors worked with permutation tableaux.

\subsubsection{Packed tableaux}

As already noticed, the arrows of a tableau are in bijection with its non free rows and columns. This implies immediately that the tableaux in $\mA_{i,j}(n)$ have exactly $n-i-j$ arrows. The maximum of $n$ arrows, i.e. the  case $i=j=0$, cannot actually be attained if $n>0$: indeed, if for instance a tableau has no free row, then the leftmost column of this tableau cannot contain any up arrow and thus is free. A total of $n-1$ arrows in a tableau can actually be reached, and in fact constitutes a fundamental class of tableaux as we will see:

\begin{defi}[Packed tableau]
\label{def:packed}
A \emph{packed tableau} of length $n>0$ is a tableau with $n-1$ arrows. Equivalently, it is a member of either $\mA_{0,1}(n)$ or $\mA_{1,0}(n)$.
\end{defi}

In fact, the following proposition shows that the unique free column of a tableau in $\mA_{0,1}(n)$ is necessarily the leftmost one:

\begin{prop}\label{prop:packed}
If $n>1$ and $T$ is a tableau in $\mA_{0,1}(n)$, the top left cell $c$ of $T$ contains a left arrow.
\end{prop}

 \dem First note that the tableau $T$ has at least one cell (so that $c$ is well defined) since tableaux with no cells have at least one free row or two free columns (note that we chose $n>1$). $T$ has no free rows, so there is a left arrow in the top row in particular. The column where this arrow lies is necessary free, because any up arrow in it would violate the alternative tableau property. But the leftmost column of $T$ is also free, because the presence of any up arrow in it would force the row where this arrow lies to be free, which is excluded. As $T$ has just one free column, this implies that $c$ contains indeed a left arrow.
 \findem

By transposition, the top left cell of tableaux in $\mA_{1,0}(n), n>1$ is filled by an up arrow. But there is a simpler way to go bijectively from $\mA_{0,1}(n)$ to $\mA_{1,0}(n)$ when $n>1$: simply {\em change the filling of the cell in the top left corner from $\leftarrow$ to $\uparrow$}. This also explains why we decided to call these two sets of tableaux with the same name: up to this arrow and the case $n=1$, they are identical.

\subsubsection{Cutting rows and columns}
\label{subsub:cut}

\begin{defi}[$cut_c$ and $cut_r$] For a nonempty alternative tableau with at least one column and no empty rows, we let $cut_c$ be the operation of deleting its leftmost column (so that all row lengths decrease by one); we define $cut_r$ similarly for deleting the topmost row.
\end{defi}

\begin{figure}[!ht]
 \centering
 \includegraphics[height=4cm]{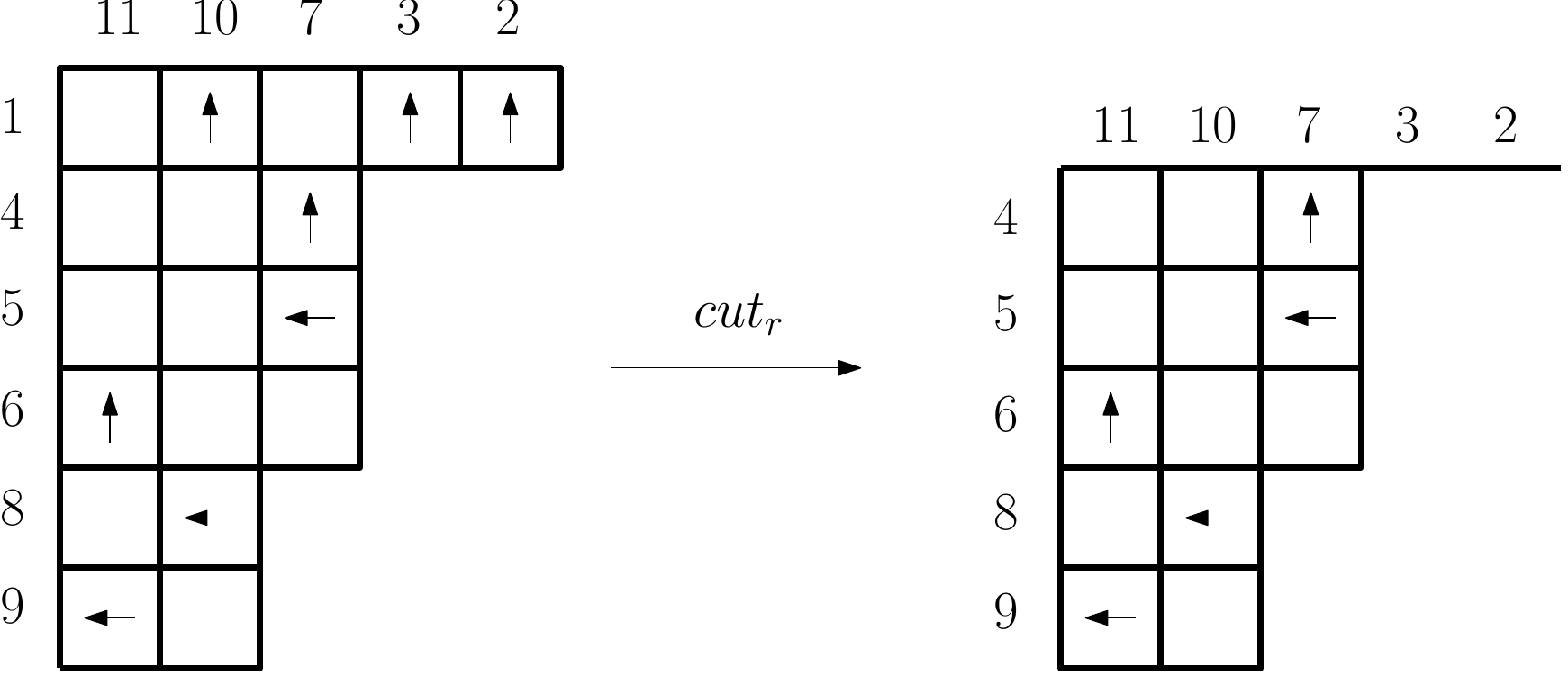}
 \caption{Cutting the first row of a tableau.
 \label{fig:cutr}}
\end{figure}

 When the tableau from which we start is labeled, we obtain naturally a labeled tableau as a result by simply keeping the labels of the remaining rows and columns in each case. This is illustrated on Figure~\ref{fig:cutr}.

  Given a tableau $T$ in $\mA$, add a new top row, and fill by up arrows all cells in this row that lie above the free columns of $T$; this construction will be called $block_c$ to reflect the fact that no free column remains after its application. We define $block_r$ symmetrically.
  Then we have the following properties:

\begin{prop}
\label{prop:cut}
For all $i,n\geq 0$, the operation $cut_r$ is a bijection between
$\mA_{i+1,0}(n+1)$ and $\mA_{i,*}(n)$; its inverse is $block_c$.
For all $j,n\geq 0$, the operation $cut_c$ is a bijection between
$\mA_{0,j+1}(n+1)$ and $\mA_{*,j}(n)$; its inverse is $block_r$.
\end{prop}

\dem
We just prove the claim concerning $cut_r$, the one for $cut_c$ being equivalent after transposition. Tableaux in $\mA_{i+1,0}(n+1)$ have no empty columns (because such columns are free), and their first row is free, since a left arrow in a cell from this row would force the corresponding column to be free. This shows that the restrictions of $cut_r$ and $block_c$ are well defined. It is immediate that $cut_r\circ block_c$ is the identity on $\mA_{i,*}(n)$. Now given a tableau $T$ in $\mA_{i+1,0}(n+1)$, we noticed that it has no left arrows in its top row, and that the up arrows in this row  occur exactly in columns that are free in $cut_r(T)$. This implies that $block_c \circ cut_r$ is the identity on $\mA_{i+1,0}(n+1)$, and the proposition is proved.
\findem

 If we do the union of the sets above for all $i$ and for all $j$ respectively, we get bijections between $\mA_{*,0}(n+1)$ and $\mA(n)$, and between $\mA_{0,*}(n+1)$ and $\mA(n)$. The special cases $i=0$ and $j=0$ are also of interest, and we get the following corollary:

\begin{cor}
\label{cor:relations}
For all $n\geq 0$, the operation $cut_r$ is a bijection between
$\mA_{*,0}(n+1)$ and $\mA(n)$, and between $\mA_{1,0}(n+1)$ and $\mA_{0,*}(n)$; $cut_c$ is a bijection between
$\mA_{0,*}(n+1)$ and $\mA(n)$, and between $\mA_{0,1}(n+1)$ and $\mA_{*,0}(n)$. Therefore for all $n\geq 0$, we have \begin{equation}\label{eq:conseqcut}
A(n)=A_{0,*}(n+1)=A_{*,0}(n+1)=A_{0,1}(n+2)=A_{1,0}(n+2)
\end{equation}
\end{cor}

\subsection{Splitting a tableau}

Now we exhibit a more complicated decomposition, which can be traced back to the last part of Burstein's work \cite{Bu}. Nevertheless in his work Burstein did not exhibit a complete recursive decomposition, mainly because he was working with permutation tableaux which are less easy to manipulate than alternative tableaux.\medskip

Let $T$ be an alternative tableau of size $n$, labeled by $L$, and $i_0$ be the label of one of its free rows (we suppose there is such a row). We compute iteratively a set of labels $T(i_0)$ in the following way: first we set $X:=\{i_0\}$. Then we add to $X$ all columns $j$ such that there is an arrow in the cell $(i_0,j)$. Afterwards, we add to $X$ all rows $i$ such that there is an arrow in $(i,j)$ for one of the columns $j$ in $X$ added at a previous stage. And so on, we keep adding row and column labels alternatively until there are no new rows or columns to add. The procedure is finite since the set $L$ is finite, and in the end we set $T(i_0):=X$. \smallskip

\noindent{\bf Example:} consider the free row labeled $4$ in the alternative tableau $T_0$ on the right of Figure~\ref{fig:permtabaltab}; so $T_0(4)$ contains $4$. Now the cell $(4,9)$ is the only cell on row $4$ containing an up arrow, so we add $9$ to the set $T_0(4)$. In this column $9$ there are two left arrow in cells $(6,9)$ and $(7,9)$, so $T_0(4)$ contains also the row labels $6$ and $7$. There is no other arrow on row $7$, and there is one on row $6$ in the cell $(6,8)$, so $8$ also belongs to $T_0(4)$. Since there is no other arrow in column $7$, we have finally $T_0(4)=\{4,6,7,8,9\}$.\smallskip

Another equivalent characterization of $T(i_0)$ is the following, which we take as a definition:

\begin{defi}[the label subset $T(i_0)$]
  Given a tableau $T$ labeled by $L$ and a free row or column $i_0\in L$, $T(i_0)$ is the smallest set $X\subseteq L$ (wrt. inclusion) which contains $i_0$, and is such that, for every cell $(i,j)$ filled by an arrow, $i$ belongs to $X$ if and only if $j$ belongs to $X$.
\end{defi}

If $Free(T)$ stands for the set of labels of free rows and columns of $T$, then we have a collection of subsets of $L$ given by $\{T(k), k\in Free(T)\}$.

\begin{lemma}
\label{lem:packed1}
Let $T\in \mA$, and $k\in Free(T)$. Then the elements of $T(k)$ other than $k$ label non free rows and columns of $T$.
\end{lemma}

\dem By the iterative definition of $T(k)$, a row label $i\neq k$ belongs to $T(k)$ if there exists a column label $j\in T(k)$ and a left arrow in the cell $(i,j)$. In particular, row $i$ is not free in $T$. The proof is similar for column labels.
\findem

\begin{lemma}
\label{lem:partition}
Let $T$ be a tableau labeled by $L$. Then the sets $T(k), k\in Free(T)$ form a partition of $L$.
\end{lemma}

\dem Suppose first that there exists an integer $p$ in $L$ that does not belong to any subset $T(k)$. We assume without loss of generality that $p$ is the label of a row, and choose the minimal such $p$. First, $p$ cannot label a free row (since it would belong to $T(p)$) so there exists $j>p$ such that $(p,j)$ contains a left arrow $\leftarrow$. Now $j$ is not the label of a free column, since otherwise $p$ would belong to $T(j)$; so there exists a row label $p'$ such that $(p',j)$ contains an up arrow $\uparrow$. We have $p'<p$ because otherwise the up arrow in $(p,j)$ would point towards the up arrow in $(p',j)$ . But then $p'$ belongs to a set $T(k)$ by minimality of $p$, which entails that $j$ and  $p$ also belong to this set, which contradicts the hypothesis that $p$ belongs to no such set.

We have thus shown that the sets $T(k), k\in Free(T)$ cover $L$, we now have to prove that they are disjoint. Let $p$ belong to a set $T(k)$; we will show that we can uniquely determine $k$ from $p$. If $p$ is free, then $k=p$ because there is only one free row or column in $T(k)$ by Lemma~\ref{lem:packed1}. Now suppose $p$ is not free, and let us assume that $p$ labels a row: there exists a (necessarily unique) column $j\in T(k)$ with $(p,j)$ containing a left arrow. Now if $j$ is free, we know that $k=j$ and we are done. Otherwise we have an up arrow in $(i,j)$ for a unique $i\in T(k)$. If $i$ is free, then $k=i$ and we are done, otherwise we continue this process, and stop until we hit upon a free row or column, and we know this is the index $k$. To conclude, we just need to be sure that the process will end: this is indeed the case because the row labels that we encounter are strictly decreasing (and the column labels strictly increasing).

\findem

Given a tableau $T$ with label set $L$, and any subset $A\subseteq L$, one can form a new tableau by selecting in $T$ only the rows and columns with labels in $A$:

\begin{defi}[$T{[}A{]}$ and $T{[}k{]}$]
Let $T$ be a tableau labeled by $L$, and $A\subseteq L$. The tableau $T[A]$ is defined as the tableau labeled by the subset $A$, where $l\in A$ labels a row (respectively a column) in $T[A]$ if it labels a row (resp. a column) in $T$, and such that the cell $(i,j)\in T[A]$ has the same filling as the cell $(i,j)$ in $T$. We write $T[k]:=T[T(k)]$ for simplicity if $k$ is a free row or a free column.
\end{defi}

 Then from Lemma~\ref{lem:packed1} we deduce immediately:

\begin{prop}
\label{prop:freeTk}
Let $T$ be a labeled tableau and $k\in Free(T)$. The tableau $T[k]$ is a packed tableau in which $k$ labels the only free row or column.
\end{prop}


\subsection{Merging tableaux}

 Since we described a way to split a tableau into smaller tableaux, it is natural to try to reconstruct the original tableau, so we need to define a way to merge tableaux together.

 \begin{defi}[the function $\merge$] Let $T$ and $T'$ be two alternative tableaux labeled on \emph{disjoint} integer sets $L$ and $L'$. Then $T''=\merge(T,T')$ is a labeled tableau defined as follows: its label set is $L''=L\cup L'$, where $k\in L''$ labels a row in $T''$ if and only if it labels a row in either $T$ or $T'$. Then the cell $(i,j)\in T''$ is filled with a left arrow if one of the following two cases occur: either $i,j\in L$ and $(i,j)$ is a left arrow in $T$, or  $i,j\in L'$  and $(i,j)$ is a left arrow in $T'$. Up arrows in $T''$ are defined similarly, and the other cells are left empty.
 \end{defi}

\begin{figure}[!ht]
 \centering
 \includegraphics[width=\textwidth]{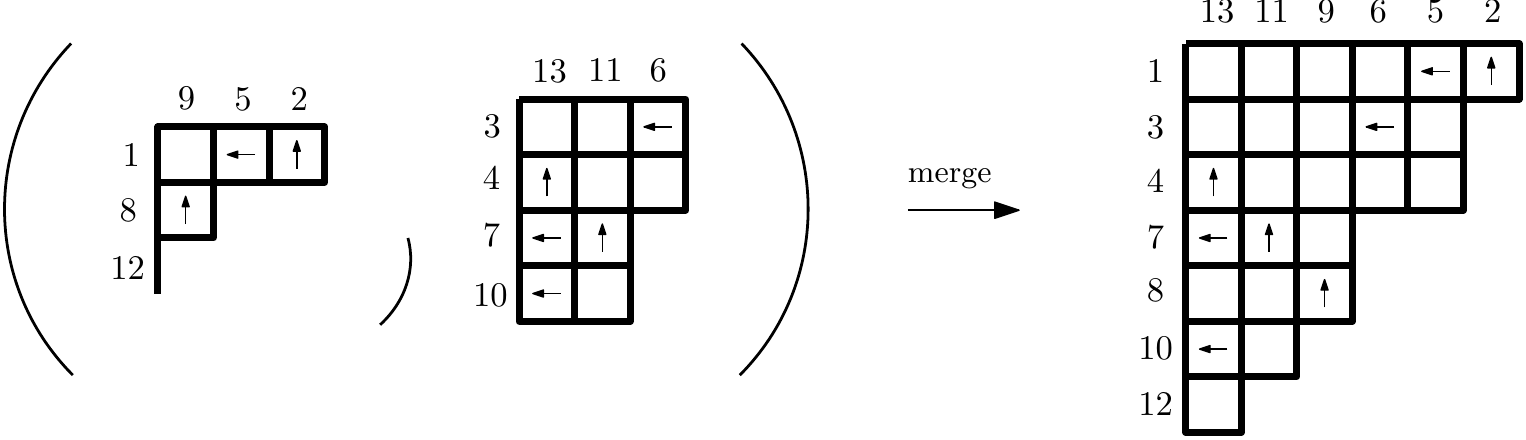}
 \caption{Merging of two tableaux.
 \label{fig:merge}}
 \end{figure}

   Note that the empty cells of $\merge(T,T')$ correspond either to empty cells in $T$ or $T'$, or to cells $(i,j)$ for which one of $i,j$ belongs to $L$ and the other to $L'$. An example of merging is given on Figure~\ref{fig:merge}.

We make a slight abuse of notation in writing $\merge(T,T')$, since the operation of merging depends crucially on the labels and not merely on the tableaux. This will not cause any problem in the rest of the paper, since we will always use it when the labels of the tableaux are clear from the context. We now record some immediate properties of the merging procedure in the following proposition:

\begin{prop}

\label{prop:merge}
Given $T,T'$ as above, then $T''=\merge(T,T')$ is an alternative tableau.
If $k\in Free(T'')$, then either $k\in Free(T)$ and $T''[k]=T[k]$, or
$k\in Free(T')$ and  $T''[k]=T'[k]$.

Moreover, the application $\merge$ is \textit{symmetric}, i.e. $\merge(T,T')=\merge(T',T)$; it is also \textit{associative}, in the sense that if $T_1,T_2,T_3$ are tableaux labeled by pairwise disjoint sets, then $\merge(T_1,\merge(T_2,T_3))=\merge(\merge(T_1,T_2),T_3)$.

\end{prop}

The last two properties allow us to extend the domain of definition of $\merge$: given a finite collection of tableaux $C=(T_i)_{i\in I}$ with pairwise disjoint label sets, we can merge all tableaux in $C$ by defining $$\merge(C)=\merge(T_{i_1},\merge(T_{i_2},\ldots, \merge(T_{i_{t-1}},T_{i_t}))),$$ where $i_1,\ldots,i_t$ is any ordering of the index set $I$; this is well defined thanks to the properties of symmetry and associativity.

\subsection{Decomposition of tableaux}

The following theorem is the main structural result of this paper; together with Corollary~\ref{cor:relations}, it describes a recursive structure that completely characterizes alternative tableaux. This will be applied in the remaining sections, first to easily obtain old and new enumerative results, and then to give bijections between alternative tableaux, certain classes of trees, and permutations of integers.

\begin{thm}
\label{th:decomp}
Let $i,j$ be nonnegative integers, and $L$ be a label set. The function $\spli:T\mapsto \{T[k],~k\in Free(T)\}$ is a bijection between:
\begin{enumerate}
\item Tableaux in $\mA_{i,j}$ labeled by $L$, and
\item Sets of $i+j$ packed tableaux, with $i$ of them in $\mA_{1,0}$ and $j$ in $\mA_{0,1}$, all labeled in such a way that their $i+j$ label sets form a partition of $L$.
    \end{enumerate}
The inverse bijection is the operation $\merge$.
\end{thm}

\dem First, the fact that $\spli$ is well defined is a consequence of Lemma~\ref{lem:partition} and Proposition~\ref{prop:freeTk} . We have also $\spli\circ \merge$ is equal to the identity function thanks to Proposition~\ref{prop:merge}. What remains to be proven is that $\merge\circ \spli$ is the identity on $\mA(n)$: that is, we need to show that given a tableau $T$, merging the labeled tableaux $\{T[k],~k\in Free(T)\}$ gives back the tableau $T$. Let us then denote by $T'$  the tableau  $\merge\circ \spli (T)=\merge((T[k])_{k\in Free(T)})$, and show that we have $T'=T$. We note immediately that the (labeled) shapes of $T$ and $T'$ coincide, so we have to show that the contents of all cells are identical. Let then $c$ (respectively $c'$) be the content of a cell $(i,j)$ in $T$ (\emph{resp.} in $T'$). If $i$ and $j$ are labels of the same tableau $T[k]$ (for a certain $k$), then $c$ is the content of $(i,j)$ in $T[k]$ by the definition of $\spli$; but by definition of $\merge$, this is also equal to $c'$. Otherwise, $i$ and $j$  belong respectively to tableaux $T[k]$ and $T[k']$ with $k\neq k'$, and in this case $c$ is necessarily empty by Lemma~\ref{lem:partition}; and by the definition of $\merge$ again, $c'$ is also empty. Thus $T=T'$ and the result is proved.

\findem

\begin{figure}[!ht]
\centering
 \includegraphics[width=\textwidth]{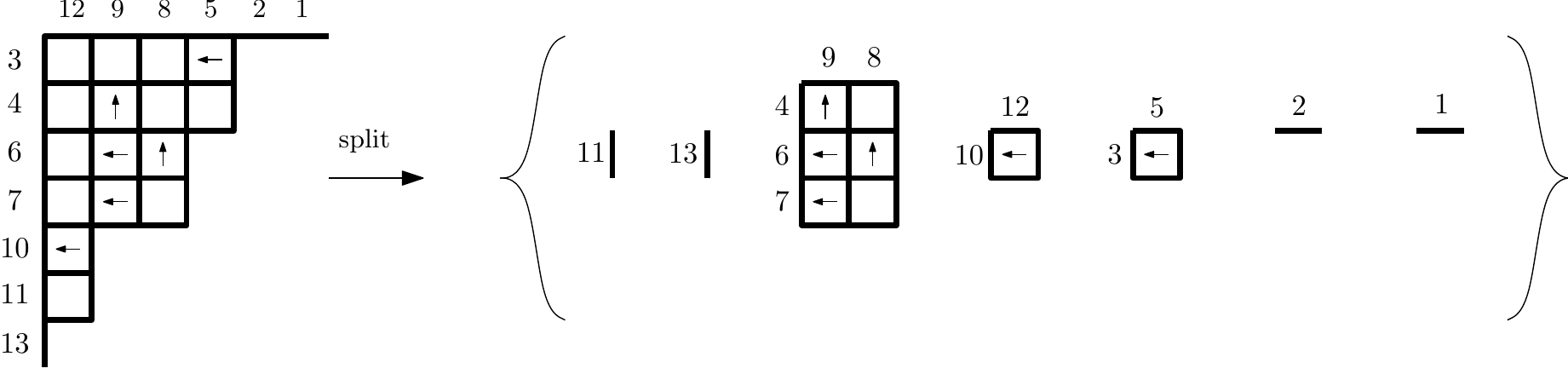}
 \caption{The decomposition $\spli$.
 \label{fig:decomptab}}
 \end{figure}

An immediate corollary of the theorem is the following, which gives a different way of decomposing tableaux:

\begin{cor}
\label{cor:decompo}
There is a bijection $\operatorname{divide}$ between $(i)$ tableaux in $\mA_{i,j}(n)$ labeled by a set $L$, and $(ii)$ pairs of tableaux $(P,Q)\in\mA_{i,0}\times \mA_{0,j}$ labeled by sets $L_P$ and $L_Q$ such that $\{L_P,L_Q\}$ is a partition of $L$.
\end{cor}

\dem Let $T$ be a tableau in $\mA_{i,j}(n)$ labeled by a set $L$. First use the bijection $\spli$ of the previous theorem, and, among the tableaux obtained, separate the ones in $\mA_{1,0}$ and the ones in $\mA_{0,1}$; merge separately each of these two collections to obtain the tableaux $P$ and $Q$ of the theorem. \findem

There is a more direct way to obtain the same bijection: consider the subsets of labels $A=\cup_k T(k)$ and $B=\cup_l T(l)$, where $k$ (respectively $l$) goes through the labels of the free rows of $T$ (resp. the free columns). Then define simply $P:=T[A]$ and $Q:=T[B]$.

\begin{figure}[!ht]
 \centering
 \includegraphics[width=\textwidth]{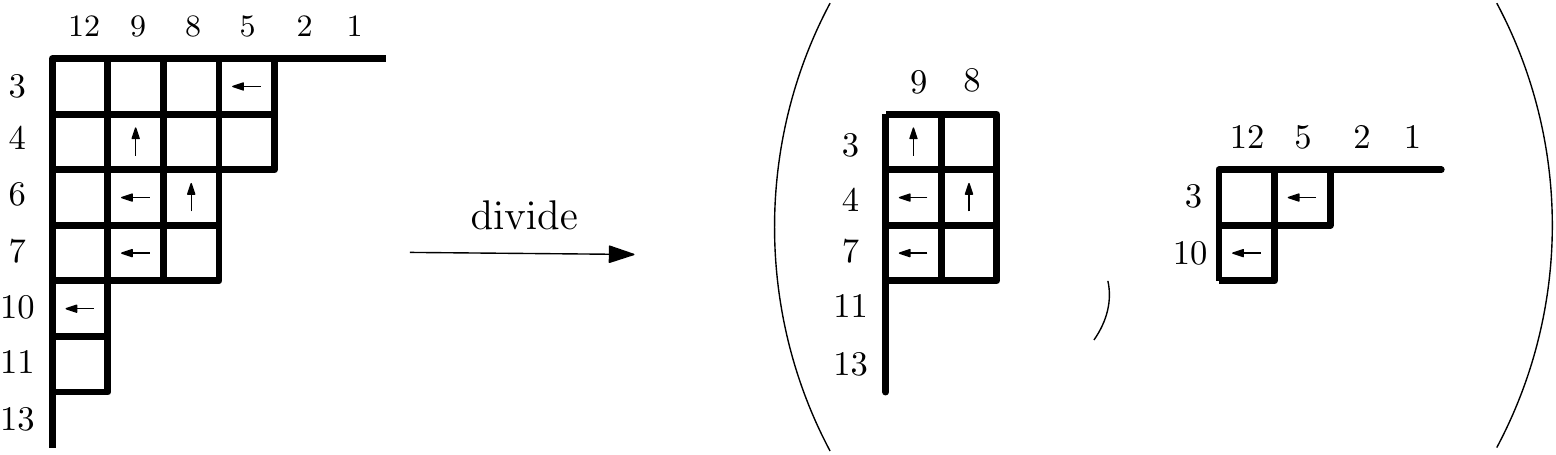}
 \caption{The tableaux $(P,Q)=\operatorname{divide}(T_0)$ for the tableau $T_0$ of Figure~\ref{fig:permtabaltab}, left.
 \label{fig:divide}}
 \end{figure}

\section{Enumeration}
\label{sect:enum}

We will show that, using the structure of alternative tableaux discovered in Section~\ref{sect:struct}, it is easy to prove various enumeration results in a simple way, starting with the plain enumeration of alternative tableaux according to their size.

\subsection{Labeled combinatorial classes}
\label{sub:combclass}
From the decompositions of Theorem~\ref{th:decomp} and Corollary~\ref{cor:decompo}, one can easily write down equations for the combinatorial class $\mA$ of alternative tableaux, in the manner of Flajolet and Sedgewick \cite{Flaj}. Indeed, Theorem~\ref{th:decomp}  says that the number of tableaux labeled on a set $L$ is the same as the number of ways to partition $L$ and then choose, for each block $b$ of this partition, a tableau labeled on $b$ belonging to either $\mA_{0,1}$ or $\mA_{1,0}$; in the language of \cite{Flaj}, this is written:
\begin{equation}
\label{eq:bla1}
\mA=SET(\mA_{0,1}+\mA_{1,0}).
\end{equation}
Similarly, a consequence of Corollary~\ref{cor:decompo} is
\begin{equation}
\label{eq:bla2}
\mA=\mA_{0,*}\star \mA_{*,0}.
\end{equation}

This means that an alternative tableau labeled on $L$ is obtained by choosing two alternative tableaux $P,Q$ in $\mA_{0,*}$ and $\mA_{*,0}$, labeled respectively by $L_P$ and $L_Q$ which are disjoint and whose union is equal to $L$. The advantage of describing our theorems in this way is that there is an automatic way to write down equations for the corresponding exponential generating functions, with the added possibility of taking into account certain parameters; this is what we will do in the rest of this Section.

As a matter of fact, the natural framework for the study of alternative tableaux is arguably the theory of {\em species on a totally ordered set}, cf. \cite[Chapter 5]{BLL}. This is not needed for the results in this work and therefore we will not develop this approach.

\subsection{The number of alternative tableaux.}
\label{sub:number}
We will give first a simple proof of the well known fact that alternative tableaux of size $n$ are enumerated by $(n+1)!$; not surprisingly, that is how the original permutation tableaux got their name. Let $A(z)$,$B(z)$ and $C(z)$ be the exponential generating functions of tableaux in $\mA$, $\mA_{0,*}$ and $\mA_{0,1}$ according to their length, that is:
$$A(z)=\sum_{n\geq 0} A(n)\frac{z^n}{n!},~
B(z)=\sum_{n\geq 0} A_{0,*}(n)\frac{z^n}{n!},\text{~and~} 
C(z)=\sum_{n\geq 0} A_{0,1}(n)\frac{z^n}{n!}.$$

On the one hand, Corollary~\ref{prop:cut} implies the following relations on generating functions:
 \begin{equation}
\label{eq:diffrel}
B'(z)=A(z)\text{~~and~~}C''(z)=A(z).
\end{equation}

 On the other hand, note that $\mA_{1,0}$ and $\mA_{*,0}$ have respectively the generating functions $C(z)$ and $B(z)$: this is immediate by transposition (cf. Proposition~\ref{prop:transp}). So we can use the combinatorial equations~\eqref{eq:bla2} and~\eqref{eq:bla1} to obtain the functional equations $A(z)=B(z)^2$ and $A(z)=\exp(2C(z))$, by an application of the principles found in \cite[Chapter II]{Flaj}.  Together with~\eqref{eq:diffrel}, we get the differential equations $$B'(z)=B(z)^2 \quad\text{and}\quad C''(z)=\exp(2C(z)).$$ With the obvious initial conditions $B(0)=1,C(0)=0, C'(0)=1$, the solutions to these are respectively $B(z)=\frac{1}{1-z}$ and $C(z)=-\log(1-z)$. Taking coefficients, we obtain $A_{0,*}(n)=n!$ and $A_{0,1}(n)=(n-1)!$, which both give us $A(n)=(n+1)!$ by Corollary~\ref{cor:relations}. To sum up we have:

\begin{prop}
\label{prop:simpleenum}
We have the following expressions:
$$
A(z)=\frac{1}{(1-z)^2},~~B(z)=\frac{1}{1-z},\text{~and~}C(z)=-\log(1-z)
$$
\end{prop}

\subsection{Refined enumeration}

In fact we can do much better by introducing some statistics. Let $A_{i,j}(n,k)$ be the number of tableaux in $\mA_{i,j}(n)$ with $k$ rows, where we allow $i=*$ or $j=*$. We define the corresponding generating functions $A_{i,j}(z,u)=\sum_{n,k\geq 0}A_{i,j}(n,k)\frac{z^n}{n!}u^k$ and $A(z,u,x,y)=\sum_{i,j\geq 0}x^iy^jA_{i,j}(z,u)$. We have then the following refined enumeration:
\begin{thm}
\label{th:enum}
\begin{equation}\label{eq:enum}
A(z,u,x,y)=\exp\left(zy(1-u)+(x+y)\ln\left(\frac{1-u}{1-u\exp(z(1-u))}\right)\right)
\end{equation}
\end{thm}

\dem By  Theorem~\ref{th:decomp}, we know that the number of free rows (respectively free columns) of a tableau is equal to the number of  tableaux in $\mA_{1,0}$ (resp. $\mA_{0,1}$) under the bijection $\spli$. We can then use Equation~\eqref{eq:bla1} and insert parameters $x$ and $y$ in it (cf. \cite[Chapter III]{Flaj}) and this gives the equation:
\begin{equation}
\label{eq:de}
A(z,u,x,y)=\exp (xA_{1,0}(z,u))\exp(yA_{0,1}(z,u))
\end{equation}

 Now we have $A_{0,1}(n,k)=A_{1,0}(n,k)$ if $n>1$, by the remark following Proposition~\ref{prop:packed}. Taking into account $n=1$, we obtain on the level of generating functions  $A_{0,1}(z,u)=A_{1,0}(z,u)+z(1-u)$; plugging into Equation~\eqref{eq:de} gives
\begin{equation}
\label{eq:de2}
A(z,u,x,y)=\exp\left((x+y)A_{1,0}(z,u)+zy(1-u)\right).
\end{equation}
 
Using the bijections $cut_r$ and $cut_c$, we get the following refinements of Corollary~\ref{cor:relations}
 $$A(n,k)=A_{0,*}(n+1,k)=A_{*,0}(n+1,k+1)=A_{1,0}(n+2,k+1),$$
for $n,k\geq 0$. This translates into the following equations for the generating functions, where all derivatives here are taken with respect to the variable $z$:
\begin{align}
A'_{0,*}(z,u)&=A(z,u);\label{eq10}\\
A_{*,0}(z,u)&=uA_{0,*}(z,u)+(1-u);\label{eq11}\\
A'_{1,0}(z,u)&=uA_{0,*}(z,u).\label{eq12}
\end{align}

Since the number of rows of $\merge(T,T')$ is the sum of the number of rows of $T$ and $T'$, we have the equation $A(z,u)=A_{0,*}(z,u)\cdot A_{0,*}(z,u)$ by Equation~\eqref{eq:bla2}. Using Equations~\eqref{eq10} and~\eqref{eq11} we get $A'_{0,*}(z,u)=A_{0,*}(z,u)\cdot (uA_{0,*}(z,u)+1-u)$. Taking into account the initial condition $A_{0,*}(0,u)=1$, this differential equation is easily solved and gives us
\begin{equation}
\label{eq:ao1}
A_{0,*}(z,u)=\frac{(1-u)}{\exp(z(u-1))-u}.
\end{equation}

Now we use Equation~\eqref{eq12}, and by immediate integration of~\eqref{eq:ao1} we obtain
$$A_{1,0}(z,u)=\ln\left(\frac{1-u}{1-u\exp(z(1-u))}\right).
$$

Now it suffices to replace $A_{1,0}(z,u)$ in~\eqref{eq:de2} and the result follows.

\findem

From this theorem, we have the following corollary, first proved in \cite{CN} by a complicated recurrence:

\begin{cor}
Define the polynomial $A_n(x,y)=\sum_{i,j}A_{i,j}(n)x^iy^j$; then
 we have the following expression:
\begin{equation}
 A_n(x,y)=\prod_{i=0}^{n-1}(x+y+i)
\end{equation}
\end{cor}

\dem It is easily seen that for $u=1$ the expression inside the logarithm in~\eqref{eq:enum} boils down to $\frac{1}{1-z}$, so

\begin{align*}
A(z,1,x,y)& =\exp(-(x+y)\log(1-z))=(1-z)^{-(x+y)}\\
 &= \sum_{n\geq 0} (x+y)(x+y+1)\cdots (x+y+n-1) \frac{z^n}{n!}.
\end{align*}

It suffices to take the coefficient of $\frac{z^n}{n!}$ on both sides to obtain the result. Note that in fact we just need Equation~\eqref{eq:de} from the proof of Theorem~\ref{th:enum}, and then the expression of $A(z,1,x,y)$ follows from  the fact that both $A_{1,0}(z)$ and $A_{0,1}(z)$ are equal to $-\log(1-z)$ by Proposition~\ref{prop:simpleenum}.

\findem

\subsection{Decorated tableaux}

In their study of the ASEP model in the case $q=1$, Corteel and Williams \cite{CW2} managed to express the stationary distribution in terms of alternative tableaux with certain weights. In particular, the so called  {\em partition function} can be expressed combinatorially. Following \cite{CW2}, let us call {\em decorated alternative tableau} an alternative tableau where each arrow can be in two states, marked and unmarked: a usual alternative tableau with $k$ arrows thus gives rise to $2^k$ different decorated alternative tableaux.

\begin{thm}[\cite{CW2}]
\label{th:nonfree}
The number of decorated alternative tableaux of length $n$ is equal to $2^n n!$.
\end{thm}

We give a simple proof of this fact based on the recursive structure of tableaux:\smallskip

\dem Let $\widetilde{A}(z)$ (respectively $\widetilde{C}(z)$) be the exponential generating function of decorated tableaux (resp.  decorated tableaux such that the underlying alternative tableau belongs to $\mA_{0,1}$). Note that, by transposition, $\widetilde{C}(z)$ can equivalently be defined by replacing $\mA_{0,1}$ by $\mA_{1,0}$. Remember that arrows correspond to non free rows and columns: so this is an additive parameter of tableaux with respect to the decomposition $\spli$. Thus from Equation~\eqref{eq:bla1} we get immediately
\begin{equation}
\label{eq:nonfree}
\widetilde{A}(z)=\exp(2\widetilde{C}(z)).
\end{equation}

But since tableaux in $A_{0,1}(n)$ (for $n\geq 1$) have exactly $n-1$ non free rows and columns, each of them gives rise to $2^{n-1}$ decorated tableaux. In terms of generating functions this means that $\widetilde{C}(z)=\frac{C(2z)}{2}$. Now we know that $C(z)=-\log(1-z)$ by Proposition~\ref{prop:simpleenum}, so after substituting in Equation~\eqref{eq:nonfree} we get $$\widetilde{A}(z)=\frac{1}{1-2z}$$ and the result follows by taking the coefficient of $z^n/n!$ on both sides.

\findem

The proof in \cite{CW2} is more involved, but has the nice feature of being bijective. It turns out that the proof above can be easily ``bijectivized'':

\begin{prop}
There is a bijection between $(i)$ decorated tableaux of length $n$, and $(ii)$ tableaux in $A_{0,*}(n)$ such that all rows and columns can be marked.
\end{prop}

Since we will give in Section~\ref{sect:perm} a bijection between tableaux of $A_{0,*}(n)$ and permutations on $n$ elements, this will indeed give a fully bijective proof of Theorem~\ref{th:nonfree}.
\smallskip

\dem
Let $T$ be a tableau of length $n$, with standard labeling, and let $P,Q$ be the tableaux respectively in $A_{*,0}$ and $A_{0,*}$ obtained by the procedure $divide$ of Corollary~\ref{cor:decompo}, together with their label sets $L_P$ and $L_Q$: we also naturally let rows and columns of $P$ and $Q$ be marked whenever they were originally marked in $T$. Now define a marked labeled tableau $P'$ as follows: the underlying tableau is $tr(P)$ and the labels are given by $L_P$. For the marks, note that transposition exchanges (free) rows and (free) columns. We keep the marks of $P$ in $P'$ for all non free rows and columns. Now $P$ has no free columns (so that $P'$ has no free rows), and all its free rows are unmarked by the definition of decorated tableaux: the corresponding free columns in $P'$ are defined to be all marked.

$P'$ and $Q$ are two marked, labeled tableaux in $A_{0,*}$, therefore $T':=\merge(P',Q)$ is a marked, labeled tableau in $A_{0,*}(n)$, and we claim that $T\mapsto T'$ is the desired bijection. Indeed, let us describe the inverse bijection. Given $U$ a marked tableau in $A_{0,*}(n)$, let $(j_1,\ldots,j_k)$ be the labels of the marked free columns, and $(l_1,\ldots, l_t)$ the labels of the unmarked free columns. Define then $R$ (respectively $S$) as the tableau $T[X]$ where $X$ is the subset of labels $\cup T(j_i)$ (resp. $\cup T(l_i)$); these are both labeled, marked tableaux in $A_{0,*}$. Now transpose the tableau $R$, keeping all marks except the ones corresponding to the original labels $(j_1,\ldots,j_k)$ which are deleted. Merge the resulting tableau $R'$ with $S$, and let $U'$ be the resulting labeled, marked tableau: it is clear that it has no marks on free rows and columns, and $U'$ is thus a decorated tableau. It is then easy to see that $U\mapsto U'$ is the wanted inverse bijection.
\findem

\subsection{Symmetric tableaux}
We call a tableau {\em symmetric} if it is fixed by the operation
of transposition defined in Section~\ref{sect:first}. Clearly
symmetric tableaux have even length since they have the same number of rows and columns. We have then the following enumeration

\begin{prop}
\label{prop:symmtab}
The number of symmetric tableaux of size $2n$ is $2^nn!$.
\end{prop}

 \dem Let $T$ be a symmetric tableau of size $2n$ with standard labeling. If $k$ labels a row, then $2n+1-k$ labels a free column. In fact, even more is true: the tableau $T[2n+1-k]$ labeled by $T(2n+1-k)$ is the transpose of the tableau $T[k]$ labeled by $T(k)$, and the labels verify $T(2n+1-k)=\{2n+1-\ell, \ell \in T(k)\}$. By the bijection of Corollary~\ref{cor:decompo}, symmetric tableaux are thus in one to one correspondence with pairs of labeled tableaux $(P,Q)\in \mA_{*,0}(n)\times \mA_{0,*}(n)$, where $Q=tr(P)$ and the labels verify 
$L_Q=\{2n+1-\ell, \ell \in L_P\}$ as well as $L_Q=\{1,\ldots,2n\}-L_P$. 

\begin{figure}[!ht]
\centering
 \includegraphics[height=3.5cm]{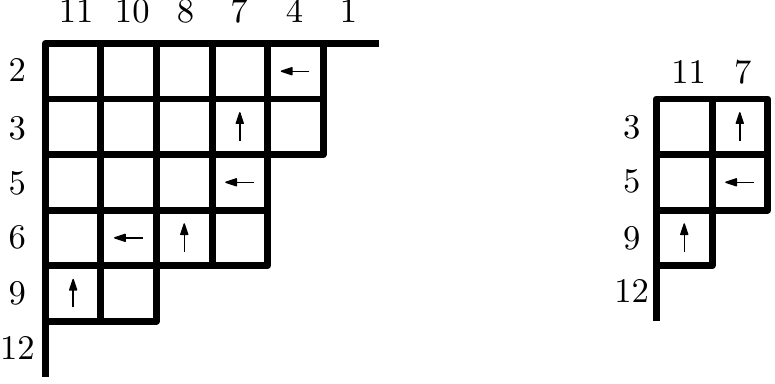}
 \caption{A symmetric tableau and its associated column packed tableau.
 \label{fig:symm}}
\end{figure}

 Thus all symmetric tableaux are obtained in the following manner: pick an alternative tableau $U$ in $\mA_{*,0}(n)$, and for each pair $\{i,2n+1-i\}, i=1\ldots n$ pick one of the two integers; let $X$ be the set of the chosen integers and $Y$ the complement of $X$ in $\{1,\ldots,2n\}$. Then merge $U$ labeled by $X$ and $tr(U)$ labeled by $Y$: this is a symmetric tableau. Since $\mA_{*,0}(n)$ has $n!$ elements and there are clearly $2^n$ choices for the labels $X$, the result follows. We illustrate the correspondence $T\mapsto P$ on Figure~\ref{fig:symm}.\findem

So symmetric tableaux of size $n$ are equinumerous with \emph{signed permutations} of $\{1,\ldots,n\}$, which are permutations on $\{1,\ldots,n\}$ in which letters can be {\em barred}; we actually describe a bijection at the end of Section~\ref{sub:bijperm}.

 In the work of Lam and Williams \cite{LW}, {\em permutation tableaux for the type $B_n$} are defined, also enumerated by $2^n n!$. It is actually possible to show that their tableaux are in bijection with symmetric alternative tableaux, by adapting suitably the bijection $\alpha$ from Theorem~\ref{th:permalt} to the symmetric case.

 We note that these permutation tableaux for the type $B_n$ were defined as a certain subclass of diagrams that appeared naturally in the context of Coxeter groups of type $B_n$. It is quite surprising that when one starts with alternative tableaux (originally related to the ASEP), and then considers those that are symmetric, one obtains configurations that are in simple bijection with permutation tableaux of type $B_n$. This raises the following question: given a finite Coxeter system $(W,S)$, is there a natural way to associate to each of the elements of $W$ a certain generalized alternative tableau ?

%
%
%
%
%

%
%

\section{Alternative trees and forests}
\label{sect:alttrees}

In this section we will give bijections from alternative tableaux to various families of planes and forests, bijections which are based on the decompositions of Section~\ref{sect:struct}.
 \smallskip

All trees considered are rooted and \emph{plane}, by which we mean as usual that the children of every vertex are linearly ordered. Furthermore, we will consider \emph{labeled} trees and forests, where the labels will be pairwise different integers attached to the vertices; these integers form the \emph{label set} of the tree or forest.

 \begin{defi}[Minimal and Maximal vertices]
 Given a vertex $v$ in a labeled tree, we say that $v$ is {\em minimal} (respectively {\em maximal}) if its label is smaller (resp. bigger) than all its descendants.
\end{defi}

\subsection{Plane alternative trees and forests}
\label{sub:pat}

\begin{defi}[Plane alternative tree]
A plane alternative tree is a labeled rooted plane tree with black and white
vertices, such that:
\begin{itemize}
\item each white vertex is minimal,
its children are black and have decreasing labels from left to right;
\item each black vertex is maximal, its children are white
 and have increasing labels from left to right.
 \end{itemize}
A \emph{plane alternative forest } is a set of plane alternative trees.
 \end{defi}

\begin{figure}[!ht]
 \centering
 \includegraphics[height=4cm]{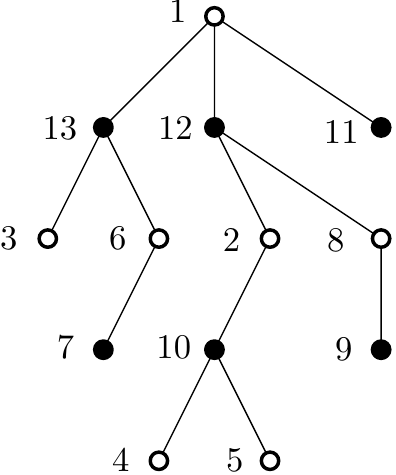}
 \caption{A plane alternative tree.
 \label{fig:planetree}}
\end{figure}

We represent a plane tree on Figure~\ref{fig:planetree}.

In this subsection we show that these trees are the natural objects encoding the recursive structure described in Section~\ref{sect:struct}.
\medskip

First we define the function $\operatorname{Tree}$ which goes from labeled \emph{packed} tableaux to alternative trees. Let $T\in \mA_{1,0}$ be labeled; then $T':=cut_r(T)$ is an element of $\mA_{0,m}$ for a certain $m$, by Corollary~\ref{cor:relations}. Let $T'_1,\ldots,T'_m$ be the labeled tableaux of $\mA_{0,1}$ given by $\spli(T')$ (cf. Theorem~\ref{th:decomp}), and $\ell$ be the label of the top row of $T$. Symmetrically, if $T\in \mA_{0,1}$, then we let the $T'_i$ be the tableaux of $\mA_{1,0}$ obtained by applying in succession $cut_c$ and $break$, and $\ell$ be the label of the leftmost column.

 \begin{defi}[$\operatorname{Tree}$ and $\operatorname{Forest}$]
 We define $\operatorname{Tree}(T)$ recursively to be the tree
  whose root is white (respectively black) and labeled by $\ell$, and
 whose subtrees attached to the roots trees are $\{\operatorname{Tree}(T'_i)\}_{i=1\ldots m}$, arranged from left to right in increasing (respectively decreasing) order of the labels of their roots if $T\in \mA_{1,0}$ (resp. $\in\mA_{0,1}$).

 If $T\in \mA$ is any labeled alternative tableau, and $\{T_i\}$ are the labeled packed tableaux given by $\spli(T)$ from Theorem~\ref{th:decomp}, we define $\operatorname{Forest}(T)$ as the labeled forest consisting of the trees $\operatorname{Tree}(T_i)$.
  \end{defi}

\begin{figure}[!ht]
 \centering
 \includegraphics[height=3cm]{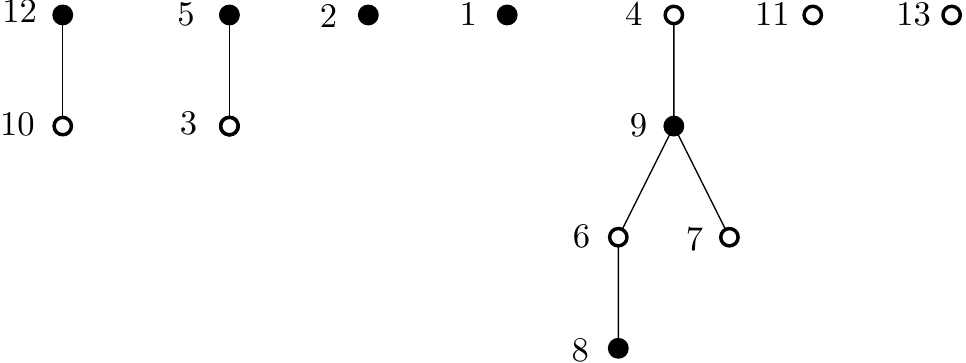}
 \caption{A plane alternative forest.
 \label{fig:arbre}}
\end{figure}

  The forest $\operatorname{Forest}(T_0)$ for the alternative tableau of Figure~\ref{fig:permtabaltab} is represented on Figure~\ref{fig:arbre}.

\begin{thm}
$\operatorname{Tree}$ is a bijection from packed labeled alternative tableaux of length $n$ to plane alternative trees with $n$ vertices (with the same label set). $\operatorname{Forest}$ is a bijection from labeled alternative tableaux of length $n$ to plane alternative forests with $n$ vertices (with the same label set).
\end{thm}

\dem Note first that the claim about $\operatorname{Forest}$ is an immediate corollary of the result for $\operatorname{Tree}$, thanks to Theorem~\ref{th:decomp}.

The proof is mostly straightforward, and consists simply of noticing that the recursive structure of tableaux given by Theorem~\ref{th:decomp} and Lemma~\ref{prop:cut} is naturally encoded by alternative trees.
The only point that needs to be checked is the minimality of white vertices (the maximality of black vertices being clearly proved symmetrically): when a white vertex is added in a tree, its label $\ell$ is the first row of a tableau $T$, and the labels of its descendants are the labels of the other rows and columns of $T$, which by definition of the labeling of tableaux are indeed larger than $\ell$.

So $\operatorname{Tree}$ is well defined, and the inverse (recursive) construction is clear: given a tree $t$ with root labeled $\ell$, construct (recursively) the labeled packed tableau $\operatorname{Tree}^{-1}(t')$ for each root subtree $t'$, then merge all these tableaux to get a tableau $T$, and finish by applying $block_c$ or $block_r$ (according to the root color) to $T$, labeling the new row or column by $\ell$: the result is $\operatorname{Tree}^{-1}(t)$.

\findem

\subsection{Arc diagrams}

We now introduce {\em alternative arc diagrams}, that turn out to be a nice representation of plane alternative forests. We will call \emph{arc diagram} the data of points aligned horizontally, labeled increasingly from left to right by integers and of arcs $(i,j), i<j$ where $i,j$ are two of the labels. It is thus a particular representation of a labeled (simple, loopless) graph where the vertices are ordered according to their value. Given an arc diagram, we say that an arc $(i,j)$ is topmost on its right side if there is no arc $(k,j)$ with $k<j$, and that it is topmost on its left side if there is no arc $(j,\ell)$ with $\ell>j$.

\begin{defi}[alternative arc diagram]
\label{def:arcdiag}
Let $L$ be a label set of size at least $2$, with minimal and maximal elements $m$ and $M$ respectively. An arc diagram with points labeled by $L$ is called {\em alternative} if the following three conditions are verified:
 \begin{enumerate}
 \item at each vertex $i$, there are no two arcs $(k,i)$ and $(i,j)$ for some integers $k<i<j$.
 \item as an abstract graph, it is a tree;
  \item each arc $(i,j)\neq (m,M)$ is topmost on exactly one of its sides.
\end{enumerate}
\end{defi}

\begin{figure}[!ht]
 \centering
 \includegraphics[height=3cm]{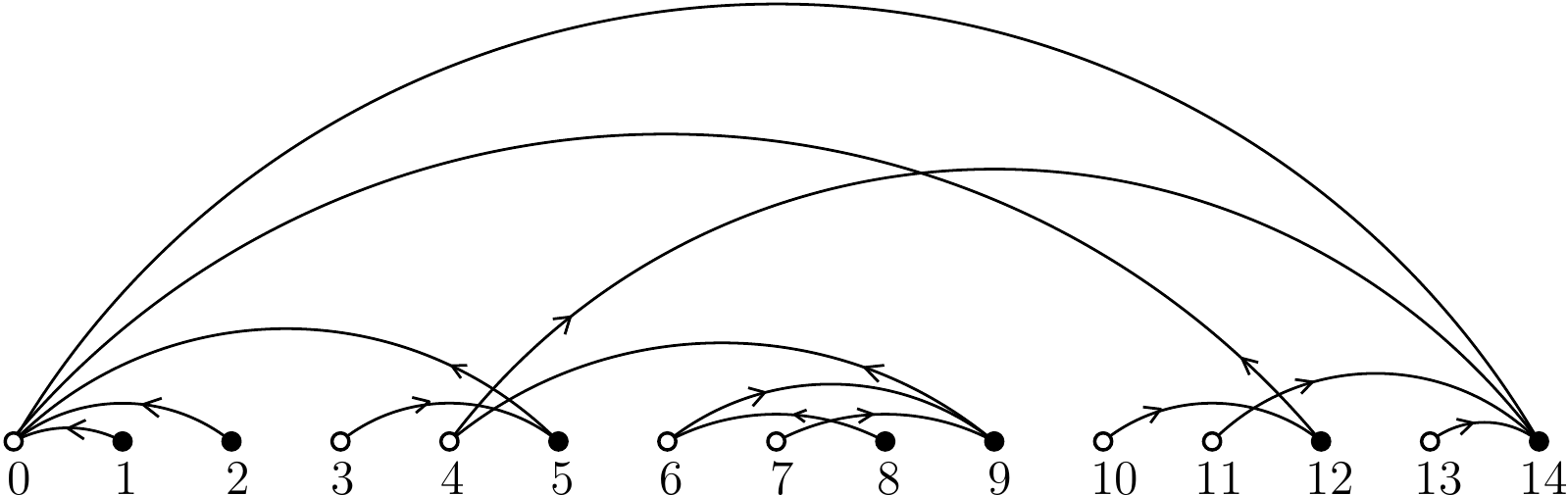}
 \caption{Alternative arc diagram.
 \label{fig:altertree}}
\end{figure}

An example is shown on Figure~\ref{fig:altertree}. Every arc has been oriented from its topmost side for clarity; moreover, a vertex $i$ is colored white when all arcs adjacent to it are of the form $(i,j)$ with $j>i$; it is colored black otherwise. Let $F$ be a plane alternative forest labeled on $\lb 1,n\rb$, and consider $n+2$ points aligned horizontally with labels $0,1,\ldots,n+1$ from left to right. Add an arc between points $i$ and $j$ for each edge $(i,j)$ of the forest, an arc $(0,b_\ell)$ for each black root $b_\ell$ and an arc $(n+1,w_k)$ for each white root $w_k$. Finally put an arc between $0$ and $n+1$, and let the resulting arc diagram be $\phi(F)$. For instance, the diagram of Figure~\ref{fig:altertree} corresponds through $\phi$ to the forest of Figure~\ref{fig:arbre}.

 \begin{prop}
The procedure $\phi$ is a bijection from alternative forests to alternative arc diagrams .
 \end{prop}

\dem
We first check the three conditions in Definition~\ref{def:arcdiag}, then show the bijectivity. So let us be given a set of arcs $\phi(F)$ on $n+2$ points coming from a forest $F$, and show that it is an arc diagram.
Condition $(1)$ is trivial for points $0$ and $n+1$; every other point $i$
is the label of a vertex in $F$. If this vertex is black, then all its descendants have smaller labels, and its father also; therefore in the diagram, all arcs go to the left of $i$. A similar proof shows that all white vertices become points from which all arcs go right. Condition $(2)$ is clear, because $F$ is a forest by hypothesis, and the arcs $(0,b_\ell)$,$(w_k,n+1)$ and $(0,n+1)$ make it into a tree.

We finally want to check condition $(3)$; let an arc $e=(i,j)\neq (0,n+1)$, $i<j$, be given. If $i=0$ or $j=n+1$ the result is immediate; now suppose $e$ is topmost in $i$; we'll show that it's not topmost in $j$, and by symmetry we'll have that if $e$ is topmost in $j$ then it's not topmost in $i$, which will conclude the proof. But $e$ topmost in $i$ means that $j$ is the father of $i$ in $F$; by minimality, the father of $j$  will necessarily be less than $i$, so that $e$ is not topmost in $j$. And if $j$ has no father in $F$, then it's a black root, thus there is an arc $(0,j)$ so that $e$ is not topmost in $j$ in this case either.

Consider now the following construction: given an arc diagram, color in white (respectively black) all points ($\neq 0,n+1$) at which arcs go to the left (resp. to the right). Then destroy all arcs $(0,j)$ and  $(j,n+1)$: the corresponding vertices $i$ and $j$ are then roots of certain trees, which form a forest. It is immediate that this is precisely the inverse of $\phi$. \findem

\subsection{Crossings in alternative arc diagrams.}

There is a very elementary way to describe the composition of $\operatorname{Forest}$ with the bijection $\phi$; let us call $Arc$ this bijection $\phi\circ Forest$ from labeled tableaux to diagrams.

 Given a tableau $T$ of length $n$ with standard labeling, and $n+2$ points labeled from $0$ to $n+1$. Then draw an arc $(i,j)$ for all cells $(i,j)$ filled with an arrow (up or left). Draw also an arc $(0,j)$ for each free column $j$, an arc $(i,n+1)$ for each free row $i$, and finally an arc $(0,n+1)$; the result is the alternative arc diagram $Arc(T)$. We have then the result

 \begin{prop}
 The construction $Arc$ is a bijection from alternative tableaux of length $n$ to alternative arc diagrams on the labels $\{0,1,\ldots,n+1\}$, and coincides with the composition $\phi\circ Forest$.
 \end{prop}

In an alternative arc diagram, we call \emph{crossing} a pairs of arcs $(i',j),(i,j')$ with $i'<i<j<j'$. Such a crossing is an \emph{out-crossing} if these arcs are topmost in $j$ and $i$ respectively. On Figure~\ref{fig:altertree}, crossings correspond to the intersection of two arcs, and out-crossings to the subset of those for which arrows are directed ``outwards'', i.e. towards $i'$ and $j'$. In this example, out-crossings occur for $(i,j)$ equal to $(4,5),(4,12),(7,8)$ and $(11,12)$. We now relate out-crossings to the free cells of alternative tableau, as defined after Definition~\ref{def:altab}, and whose importance is underlined by Proposition~\ref{prop:connection}. These free cells are also of interest in connection with permutations, see \cite{CN,SW} for instance.

\begin{prop}
\label{prop:crossings}
Let $T$ be an alternative tableau with standard labeling. A cell $(i,j)$ in $T$ is free if and only if there exists $i',j'$ such that $(i,j'),(i',j)$ is an out-crossing of $Arc(T)$; $i'$ and $j'$ are in this case unique.
\end{prop}

\dem By definition, a cell $(i,j)$ is free if the two following conditions are verified:
\begin{itemize}
\item Row $i$ is free, or there is a left arrow in a cell $(i,j')$ with $j<j'$;
\item Column $j$ is free, or there is an up arrow in a cell $(i',j)$ with
$i'<i$.
\end{itemize}

Note that the indices $j'$ and $i'$ are necessarily unique if they exist. The first condition corresponds in the arc diagram to a unique arc $(i,j')$ topmost in $i$ with $j'>i$, while the second condition corresponds to a unique arc $(i',j)$ topmost in $j$ with $i'<i$, which achieves the proof.
\findem

Note that free cells are not easily visualized when looking at plane alternative forests. As a corollary, we have the following well known enumeration, of which we give here a new simple bijective proof.

\begin{cor}\cite{CN,CW,V}
Tableaux of size $n$ with no free cells are counted by the Catalan number $C_{n+1}=\frac{1}{n+2}\binom{2n+2}{n+1}$.
\end{cor}

\dem
By Theorem~\ref{prop:crossings}, such tableaux are in bijection with alternative arc diagrams on $n+2$ points with no out-crossing. In fact, such diagrams have no crossing at all: suppose there was such a crossing $(i',j),(i,j')$ in $Arc(T)$ with $i<i'<j<j'$. Then in the tableau $T$ there are arrows in both $(i',j)$ and $(i,j')$; but this implies that the cell $(i,j)$ is free, which is absurd because this would mean that there is a out-crossing in $Arc(T)$.

 So we have to enumerate alternative arc diagrams with no crossings, and in this case Condition $(3)$ in Definition~\ref{def:arcdiag} is easily seen to be superfluous; the arc diagrams $Arc(T)$ for $T$ of size $n$ with no free cells are then identified with the well known called {\em noncrossing alternating trees} on $n+2$ points. These objects are in a simple bijection with binary trees with $n+1$ leaves, and thus are counted by the Catalan number $C_n$ : this is done in \cite{Stan2}, exercise 6.19 (p) for instance.

\findem


\subsection{Binary alternative trees}

We describe more briefly the trees that appear when one encodes the recursive structure of alternative tableaux reflected by Corollary~\ref{cor:decompo}; as can be expected, binary trees are obtained.

\begin{defi}[Binary alternative trees] A \emph{binary alternative tree} of size $n$ is a labeled binary tree with $n$ vertices such that each left child is maximal, while each right child is minimal; the root is either maximal or maximal.
 \end{defi}

We will note $\mathcal{B}_{min}$ (respectively $\mathcal{B}_{max}$) the class of binary alternative trees where the root is minimal (resp. maximal). We remark that these trees were already defined by Burstein \cite{Bu} in the context of permutation tableaux. They consist a variation of the \emph{binary increasing trees}, in which every vertex is minimal; here we distinguish left and right sons.

\begin{figure}[!ht]
 \centering
 \includegraphics[height=3cm]{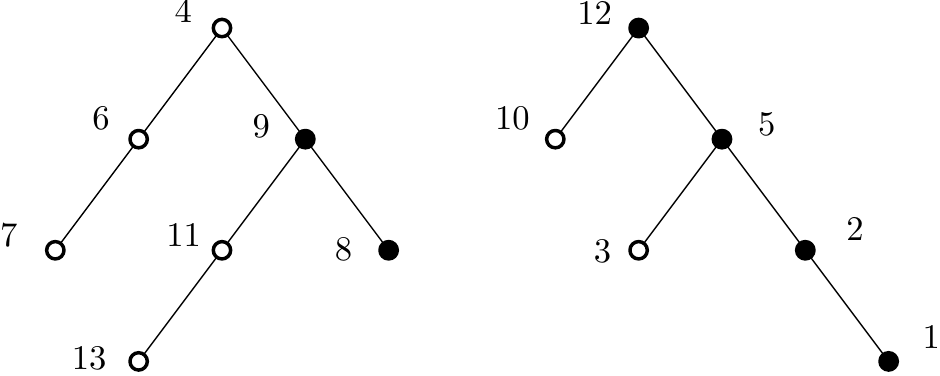}
 \caption{Binary alternative trees.
 \label{fig:bintrees}}
\end{figure}

Let $n>0$, and $T$ be a tableau in $\mA_{0,*}(n)$ labeled by $L$. If $T$ is the empty tableau, set $Bin_{min}(T)=Bin_{max}(T)=\emptyset$. Otherwise, define $m\in L$ to be the label of the first row of $T$, let $T'$ be the labeled tableau $cut_r(T)$, and finally let $P$ and $Q$ be the labeled tableaux given by $(P,Q):=divide(T')\in\mA_{*,0}\times \mA_{0,*}$ using the bijection of Corollary~\ref{cor:decompo}. By induction, we define $Bin_{min}(T)$ as the tree in $\mathcal{B}_{min}$ with a root labeled $m$ which has right subtree equal to $Bin_{min}(P)$ and left subtree equal to $Bin_{max}(Q)$.

$Bin_{max}(T)$ is defined similarly for tableaux $T$ in $\mA_{*,0}(n)$, except that $m$ is the label of the first column of $T$ and $T'=cut_c(T)$.

\begin{lemma}
$Bin_{min}$ is a bijection from $\mA_{0,*}$ to $\mathcal{B}_{min}$, and $Bin_{max}$ is a bijection from $\mA_{*,0}$ to $\mathcal{B}_{max}$.
\end{lemma}

Now let $T$ be any labeled alternative tableau, and set $$CoupleBin(T):= (Bin_{max}(P),Bin_{min}(Q))$$, in which $(P,Q)$ are the labeled tableaux given by $divide(T)$.

\begin{thm}
$CoupleBin$ is a bijection from alternative tableaux labeled
by a set $L$ to pairs of trees $(b_1,b_2)\in \mathcal{B}_{max}\times
\mathcal{B}_{min}$ with respective labels $L_1$ and $L_2$
verifying $L_1\sqcup L_2=L$.
\end{thm}

\section{Alternative tableaux and permutations}
\label{sect:perm}

In this Section we define a bijection from alternative tableaux to permutations, which relies on the representation of tableaux as trees from Section~\ref{sub:pat}. We then show that this bijection is equivalent to some other ones that already appeared in the literature.

\subsection{Some definitions}

We define a permutation as a word on the alphabet of integers with \emph{no repeated letters}. For a permutation $w=a_1a_2\cdots a_k$, we define the \emph{support} of $w$ as $supp(w):=\{a_1,\ldots,a_k\}$, i.e. the set of positive integers that appear in it; by definition  of a permutation this set has cardinal equal to $k$.


A RL-maximum (respectively a RL-minimum) in a permutation is a letter that is greater (resp. smaller) than all the letters to its left. RL stands for ``right to left'', a RL-minimum being a letter that is greater than all those seen before when one reads the word from right to left.

\begin{defi}[shifted RL-maximum]
Let $w$ be a permutation, and consider its factorization $w_1mw_2$, where $m$ is the smallest element of $supp(w)$. A  \emph{shifted RL-maximum} of $w$ is a $RL$-maximum of the permutation $w_1$.
\end{defi}

A \emph{descent} in a permutation $a_1\cdots a_k$ is a letter $a_i$ greater than $a_{i+1}$, and an \emph{ascent} is a letter smaller than the next one; by convention the last letter of a word is considered to be an ascent.


\subsection{Bijection with permutations}
\label{sub:bijperm}
We construct here a bijection $\Psi$ from plane alternative forests to permutations; composition with the function $Forest$ will give us a bijection $\Phi_N$ from tableaux of length $n$ to permutations of $\{0,\ldots,n\}$.
\smallskip

Let $T$ be a plane alternative tree; we define a permutation $\psi(T)$ recursively. If $T$ is reduced to one vertex labeled $m$, then we set $\psi(T)=m$. Otherwise, let $T_1,\ldots,T_k$ be the subtrees attached to the root (from left to right), and $m$ be the label of the root. Then the permutation attached to $T$ is the word $\psi(T):=\psi(T_1)\ldots \psi(T_k)m$. In other words, we do a {\em postorder traversal} of the tree.

\begin{lemma}
$\psi$ is a bijection between $(i)$ trees with a black root (respectively white root) labeled on $L$ and $(ii)$ permutations with support $L$ ending with the letter $\max(L)$ (resp. $\min(L)$).
\end{lemma}

\dem
The key observation is the following: if $T$ is a tree with subtrees $T_1,\ldots,T_k$ as above, then in the permutation $w=\psi(T_1)\ldots \psi(T_k)$, the last letters of the words $\psi(T_i)$ are exactly the RL-minima of $w$ (respectively the RL-maxima of $w$) if $T$ has a black root (resp. a white root). This is proved immediately by induction,  since it is a translation of the fact that black vertices are maximal, white vertices are minimal, and that the subtrees of a black vertex (respectively a white vertex) are ordered in the increasing order of their root labels (respectively the decreasing order of these labels). From this remark one can immediately defie an inverse to $\psi$.
\findem

Note that from this Lemma and the bijection $Tree$, we have that $A_{0,1}(n)=(n-1)!$ immediately, as proved in Proposition~\ref{prop:simpleenum}.
\medskip

Now let $F$ be alternative forest of size $n$ with label set $L$, composed of the trees $T_1,\ldots,T_i$ with white roots, ordered in increasing order of their roots, and $T'_1,\ldots,T'_j$ with black roots in decreasing order of their roots. Let us also fix $x<\min (L)$. Then the permutation $\Psi(F)$ is defined as the concatenation $$\Psi(F):=\psi(T'_1)\cdots \psi(T'_j)\cdot x\cdot \psi(T_1)\cdots \psi(T_i).$$

\begin{prop}
\label{prop:fortoperm}
Let $L$ be a label set, and $x<\min (L)$. Then $\Psi$ is a bijection between plane alternative forests labeled by $L$ and permutations $w$ such that $supp(w)=L\cup\{x\}$. If $\sigma=\Psi(F)$ and $i\in L$, then $i$ labels a white root (respectively a black root, a white vertex, a black vertex) of $F$ if and only if $i$ is a RL-minimum in $\sigma$ (resp. a shifted RL-maximum, an ascent, a descent).
\end{prop}

\dem The proof that $\Phi$ is bijective is essentially the same as the one for $\psi$, and the rest follows immediately from its definition.
\findem

Now we can define our bijection $\Phi_N:=\Psi\circ Forest$ from labeled alternative tableaux to permutations. Note that it requires to fix not only a labeled set $L$, but also an integer $x$ smaller than $\min(L)$. When $T\in \mA$ has the standard labeling, we will naturally take $x=0$, and with this convention we have the following theorem:

\begin{thm}[Bijection $\Phi_N$]
The bijection $\Phi_N$ is a bijective correspondence between alternative tableaux of size $n$ and permutations of $\{0,\ldots,n\}$. Furthermore, if $\sigma=\Phi_N(T)$ and $i\in \{1,\ldots,n\}$, then $i$ labels a row (respectively a column, a free row, a free column) in the standard labeling of $T$ if and only if $i$ is an ascent  (resp. a descent, a RL-minima, a shifted RL-maxima) of $\sigma$.
\end{thm}

Note that Theorem~\ref{th:enum} thus gives a refined enumeration of permutations of $\{0,\ldots,n\}$ with respect to ascents, RL-minima and shifted RL-maxima. In fact the generating function~\eqref{eq:ao1} is the generating function of Eulerian polynomials (see \cite[p.51]{Comtet}).

\medskip

 It is easy to see what \emph{symmetric tableaux} become via the bijection $\Phi_N$. Using the decomposition described in the proof of Proposition~\ref{prop:symmtab}, it is equivalent to perform the bijection $\Phi_N$ on tableaux in $\mA_{*,0}(n)$ labeled by sets $L\subseteq \{1,\ldots,2n\}$ such that $L$ contains exactly one element in each pair $\{i,2n+1-i\}$ for $i=1,\ldots,n$. The permutations obtained by $\Phi_N$  are exactly words of length $n$ labeled by such sets $L$, preceded by a $0$. Now if we delete this $0$ and replace each entry $2n+1-i$ ($i\leq n$) in this word by a barred letter $\bar{i}$, then we get a bijection with permutations where letters may be barred:

\begin{prop}
\label{prop:symmtabtoperm}
The bijection $\Phi_N$ induces a bijection between symmetric alternative tableaux of size $2n$ and signed permutations of size $n$, i.e. permutations on $\{1,\ldots,n\}$ such that each letter may be barred.
\end{prop}

This gives a bijective proof of the fact that symmetric tableaux of size $2n$ are counted by $2^nn!$.

\subsection{An ubiquitous bijection}

In this section we point out that the bijection $\Phi_N$ is identical to two  bijections that have appeared previously in the literature.

\subsubsection{Corteel and Nadeau's {\em bijection I}}

In the work of the author with Sylvie Corteel \cite{CN}, two bijections were defined between permutation tableaux and permutations; we show that the first of these bijections is identical to the bijection $\Phi_N$.

We recall this bijection $\Phi_{C}$; starting with a tableau $T$, we will define it algorithmically, by successively inserting row and column labels in a word until we reach the desired permutation. Initialize the word to the list of the labels of free rows in increasing order, preceded by $0$. Considering the columns of $T$ successively from left to right, perform the following with $j$ the current column label: if the column has no up arrow, insert $j$ to the left of $0$, while if it has an up arrow in position $(i,j)$ then insert $j$ to the left of $i$. In both cases, if $i_1,\ldots ,i_k$ are the labels of the rows containing a left arrow in column $j$, insert $i_1,\ldots ,i_k$ in increasing order to the left of $j$. When the rightmost column has been processed, we have obtain is the desired permutation $\Phi_C(T)$.

{\bf Example:} Let us apply this on the left of tableau $T_0$ of Figure~\ref{fig:permtabaltab}; the free rows are labeled by $4,11$ and $13$, so we obtain initially $(0,4,11,13)$. For column number $12$, no up arrow, a left arrow in row $10$: we get $(10,12,0,4,11,13)$. Column number $9$ has an up arrow in row $4$ and a left arrow in rows $6$ and $7$: we thus obtain $(10,12,0,6,7,9,4,11,13)$. For the remaining columns $8,5,2,1$, we obtain successively \begin{eqnarray*}
(10,12,0,8,6,7,9,4,11,13),(10,12,3,5,0,8,6,7,9,4,11,13),\\
(10,12,3,5,2, 0,8,6,7,9,4,11,13)\end{eqnarray*} and finally $$\Phi_{C}(T_0)=(10,12,3,5,2,1,0,8,6,7,9,4,11,13).$$

This is the same result as applying $\Phi_N$, and this is indeed no coincidence:

 \begin{prop}
  The bijection $\Phi_{C}$ coincides with the main bijection $\Phi_N$.
 \end{prop}

\dem
We will prove that the plane alternative forest corresponding to the permutation $\Phi_{C}(T)$ coincides with the plane alternative forest attached to an alternative tableau $T$, i.e. that we have $\operatorname{Forest}=\Psi^{(-1)}\circ \Phi_C $.

The reasoning goes by induction on the number of columns of $T$. Suppose first that $T$ has no column, and let $i_1<\ldots< i_k$ be the labels of its (necessarily free) rows. Then $\Phi_{C}(T)$ is simply the permutation $0,i_1,\ldots,i_k$, and the forest attached to this permutation is nothing else than the completely disconnected graph with $k$ white vertices labeled by $i_1,\ldots,i_k$: this is indeed the forest $\operatorname{Forest}(T)$.

Now suppose that $T$ possesses $m>0$ columns, let $j$ be the label of its rightmost column, and define $i_1<\ldots<i_k$ as the row labels of left arrows in column $j$. Let $T_1$ be the tableau obtained by suppressing this column (we keep all the labels and arrows of all other rows and columns); by induction, we know that $\sigma_1:=\Phi_{C}(T_1)$ corresponds to the forest $F_1:=\operatorname{Forest}(T_1)$. Let $\sigma:=\Phi_C(T)$ and $F:=\Psi^{-1}(\sigma)$. We distinguish two cases:

\begin{enumerate}
\item Column $j$ of $T$ has no up arrow. Then the permutation $\sigma$ is obtained by inserting $i_1\cdots i_kj$ to the left of $0$ in $\sigma_1$. The corresponding forest $F$ is obtained by adding a new black root to $F_1$ labeled $j$, and attach to it the white vertices $i_1,\ldots ,i_k$ (which were previously isolated).
\item Column $j$ of $T$ has an up arrow in row $i$. Then the permutation $\Phi_{C}(T)$ is obtained by inserting $i_1\cdots i_kj$ to the left of $i$ in the permutation $\sigma_1$. The corresponding forest $F$ is obtained by adding a new black vertex to $F_1$ labeled $j$, making it the leftmost vertex of $i$, and attach to it the white vertices $i_1,\ldots ,i_k$.
\end{enumerate}

In both cases, the forest $F$ obtained is easily seen to be precisely $\operatorname{Forest}(T)$. This proves by induction that the two functions $\operatorname{Forest}$ and $\Psi^{(-1)}\circ \Phi_C $ coincide, and thus we get indeed $\Phi_N=\Psi\circ \operatorname{Forest}=\Phi_C$.

\findem

\subsubsection{Other bijections}
 After introducing the concept of alternative tableaux in \cite{VienCamb}, Viennot defines a bijection $\Phi_V$ with permutations, which he presents under different equivalent forms. One of these consists in starting from a permutation, and the shape of a tableau (computed according to the ascents and descents of the permutation), and proceeds to fill the tableau little by little. Under this form, it is possible to show by induction it is equivalent to the bijection $\Phi_{N}$, in a similar way to what was done for$\Phi_C$ above.

At the end of Burstein's paper \cite{Bu}, a bijection is also introduced. We will not go into detail, but it is possible to see that his bijection is essentially equivalent to the other ones encountered, up some elementary transformations of permutation tableaux and of permutations.

Finally, there are two other bijections in the literature: Corteel and Nadeau's \textit{bijection II}, which is at the core of the paper \cite{CN}, and Steingr\'imsson and Williams's original bijection \cite{SW}, which is known to be equivalent to one in Postnikov's preprint \cite{Postnikov}. It would be interesting to study how these bijections are related to $\Phi_N$.

\end{document}